\documentclass[12pt,a4paper]{article}

\usepackage[utf8]{inputenc}
\usepackage[T1]{fontenc}
\usepackage{amsmath,amssymb,amsthm,amsfonts} % AMS math
\usepackage{mathtools} % extra math tools
\usepackage[left=1.2in, right=1.2in, top=1in, bottom=1in]{geometry}
\usepackage{hyperref}  % hyperlinks
\usepackage{graphicx}  % figures
\usepackage{subcaption}  % subfigure environment
\usepackage{float}
\usepackage{hyperref,url}
\usepackage{enumitem}  % better lists
\usepackage{cite}      % nice citations
\usepackage{color}     % colored text
\usepackage{bm}        % bold math symbols
\usepackage{authblk}
\usepackage[normalem]{ulem}

\newtheorem{theorem}{Theorem}[section]

\theoremstyle{definition}

\theoremstyle{remark}
\newtheorem{remark}[theorem]{Remark}

\newcommand{\eps}{\varepsilon}
\newcommand{\p}{\partial}

\newcommand{\N}{\mathbb{N}}
\newcommand{\C}{\mathbb{C}}
\newcommand{\dd}{\mathop{}\!\mathrm{d}}
\newcommand{\ee}{\mathrm{e}}
\newcommand{\ii}{\mathrm{i}}
\DeclareMathOperator{\sinc}{sinc}

\title{Numerical reconstruction of Schr\"odinger equations with quadratic nonlinearities}
\author[1]{Khaoula El Maddah}
\author[2]{Matti Lassas}
\author[1]{Teemu Tyni}
\affil[1]{Research Unit of Applied and Computational Mathematics,
University of Oulu, Finland}
\affil[2]{Department of Mathematics and Statistics, University of Helsinki, Finland}

\date{}

\begin{document}

\maketitle

\begin{abstract}

We introduce a numerical framework for reconstructing the potential in two dimensional
semilinear elliptic PDEs with power type nonlinearities from the nonlinear
Dirichlet to Neumann map. By applying higher order linearization method, we compute the Fourier
data of the unknown potential and then invert it to recover \(q\). Numerical experiments
show accurate reconstructions for both smooth and discontinuous test cases.

\end{abstract}

%\tableofcontents

\section{Introduction}

In this paper, we study the semilinear elliptic model
\begin{equation}\label{eq:nonlinear}
-\Delta u + q\,u^{p} = 0 \quad \text{in }\Omega,\qquad
u\big|_{\partial\Omega} = f,
\end{equation}
with \(p\in\mathbb{N}_{\ge 2}\) on a bounded domain \(\Omega\subset\mathbb{R}^n\) and Dirichlet data \(f\). Our objective is to recover the unknown potential \(q=q(x)\) from boundary measurements associated with \eqref{eq:nonlinear}. We study this numerically by using the higher order linearization framework of \cite{lassas2021inverse} for Calder\'on type semilinear equations with power nonlinearities. The key idea is to exploit higher derivatives
of the nonlinear Dirichlet-to-Neumann (DN) map to obtain auxiliary linearized equations,
thereby enabling recovery in settings where the corresponding linear inverse
problem is difficult to resolve. In particular, \cite{lassas2021inverse} shows that,
given the nonlinear DN map, one can simultaneously determine a potential and the
underlying conformal manifold in two dimensions, and recover a potential on
transversally anisotropic manifolds in dimensions \(n\ge 3\). In the Euclidean case,
the method yields a notably simple solution to the Calder\'on problem for certain
semilinear equations, without relying on complex geometrical optics (CGO) constructions.

We numerically implement the higher order linearization method for elliptic
equations with power type nonlinearities based on the method detailed in \cite{lassas2021inverse,liimatainen2022inverse}. The equation
\eqref{eq:nonlinear} is well-posed for boundary values \(f\in C^{2,\alpha}(\p\Omega)\), and the boundary measurements
are encoded by the nonlinear (DN) map 
\[
\Lambda_q:\ f \longmapsto \partial_\nu u\big|_{\partial\Omega},
\]
which assigns to each boundary value its corresponding normal flux. By evaluating higher Fr\'echet derivatives of
\(\Lambda_q\), we can obtain an Alessandrini-type integral identity
\[
p! \int_\Omega q\, v_0 v^p \, \dd x 
= \frac{d^p}{d\eps^p} \int_{\p\Omega} v_0 \, \Lambda_q(\eps f) \, \dd S.
\]
Particularly, it suffices to use measurements of the nonlinear functional $\psi: C^{2,\alpha}(\p\Omega)\to \C$ given by
\[
\psi:f\mapsto \int_{\p\Omega} \Lambda_q(f)\dd S_x,
\]
rather than the full DN-map.
By suitable choice of boundary values $f$, this integral identity produces the Fourier transform of the unknown potential
\(q\). Regularized Fourier inversion is then implemented to recover \(q\) from noisy and band-limited data as a solution to a regularized least squares problem (Tikhonov or total variation, depending on the expected regularity of $q$).
The methodology applies to any integer power nonlinearity $p\in\mathbb{N}_{\ge 2}$, but
we concentrate on the case \(p=2\) for simplicity.

In this work we present, to the best of our knowledge, the first numerical implementation of the higher order linearization method for elliptic nonlinear PDEs. For the forward problem, we propose a discretization scheme for the nonlinear equation \eqref{eq:nonlinear}. For the inverse problem, we design and analyze regularization strategies that combine Savitzky-Golay smoothing for differentiating the noisy data with an isotropic Fourier inversion based on Tikhonov or total variation regularization. The overall reconstruction pipeline is validated through numerical experiments on a range of test configurations, demonstrating the robustness of the method with respect to noise, discretization parameters, and model variations.

For the classical (linear) Schr\"odinger Calder\'on problem, several results
establish uniqueness and reconstruction under progressively weaker regularity. In two
dimensions, Nachman \cite{nachman1996global} gave a constructive solution for
conductivities (potentials) in \(W^{2,p}\) using \(\bar\partial\) (D-bar) techniques, giving an
explicit reconstruction formula. Nachman also addressed anisotropic formulation in
domains of \(\mathbb{R}^n\), and solve it for \(C^3\) matrix valued potentials. The planar
result was later sharpened by Astala-P\"aiv\"arinta \cite{astala2006calderon}, who proved
global uniqueness for bounded conductivities, relying on CGO
methods and removing the earlier smoothness assumptions. For higher-dimensional results (\(n\ge 3\)), we refer to Sylvester-Uhlmann \cite{sylvester1987global} and, more recently, Carstea \cite{carstea2024calder}. Other results for Calderón problem on Lipschitz domains, under the assumption that the conductivity belongs to $W^{1,\infty}(\bar\Omega)$, were established by Tataru and Haberman~\cite{haberman2013uniqueness}. These results were improved in~\cite{caro2016global, caro2025reconstruction}, where the
regularity assumptions were weakened to require only the boundedness and the Lipschitz continuity of the conductivity and where explicit reconstruction formulas were also provided.

Over the past decade, nonlinear inverse problems have progressed rapidly \cite{Uhlmann2009}, driven by the
realization that nonlinear interactions can be exploited as information rather than an obstacle. In the parabolic setting, Isakov \cite{Isa93} showed that the first linearization
of the nonlinear Dirichlet to Neumann (DN) map coincides with the DN map of a suitable linear equation,
thereby enabling the use of linear inverse theory. For the semilinear Schr\"odinger equation
\(\Delta u + a(x,u)=0\), recovery of the nonlinearity was established in 2D \cite{IN95,Sun10} and in
higher dimensions \cite{IS94}. More recently, Johansson–Nurminen–Salo \cite{johansson2025inverse,johansson2023inverse}
proved that, in a neighborhood of a fixed background solution, the entire nonlinearity \(a(x,u)\)
is determined from boundary data, up to the natural gauge. The key method often termed higher order linearization has since become a general tool, with parallel developments for
nonlinear elliptic models (including power type semilinear equations \cite{liimatainen2022inverse,gkikas2025semilinear,kian2024determination}),
quasilinear equations \cite{Sun96,SU97,KN02}, the degenerate \(p\)–Laplace equation \cite{SZ12,BHKS18}, and
and fractional semilinear Schr\"odinger problems \cite{LL19}.

Similar equations of the form~\eqref{eq:nonlinear}, often in more complex settings, arise in a variety of applications
(e.g., steady states of reaction diffusion systems with power type kinetics \cite{Smoller1994,Pao1992},
population dynamics with nonlinear birth/death terms \cite{Murray2002,OkuboLevin2001} and combustion and Emden–Fowler
type models \cite{FrankKamenetskii1969,Chandrasekhar1939,GidasNiNirenberg1979}), where \(q(x)\)
encodes spatial heterogeneity and \(p>1\) the order of the nonlinearity. Despite substantial theoretical progress, numerical implementations for nonlinear
elliptic inverse problems remain comparatively limited; our work takes a step in that direction.

The paper is organized as follows. In Section~\ref{section:HOL} we explicit the higher order linearization
method for the semilinear elliptic equation \eqref{eq:nonlinear}, as explained in \cite{lassas2021inverse}. Section \ref{section:newton} details the Newton framework used to solve the forward numerical problem using FEM. Section~\ref{section:rec method of q} focuses on the reconstruction of the
potential: We first recover Fourier data by applying the higher order linearization method to the
quadratic case \(p=2\) (Section~\ref{section:HOL}), then we specify the choice of the boundary conditions (Section~\ref{boundary choice}); next, we derive a matrix formulation
of the problem (Section~\ref{subsection: E}); finally, we invert the Fourier transform and
reconstruct \(q\) from noisy measurements (Section~\ref{subsection: q_rec}). Section~\ref{section: Numerical_examples} presents
numerical experiments. The appendices contain supporting material:  Appendix~\ref{app: Newton justified} justifies the convergence of Newton
method; and Appendix~\ref{Savitzky-Golay} discusses regularized higher order differentiation using Savitzky-Golay filter via local polynomial regression.

\section{Higher order linearization method for the semilinear elliptic equation}\label{section:HOL}
Let $\Omega \subset \mathbb{R}^n$ be a bounded domain with smooth boundary $\p\Omega$.  
We consider the semilinear elliptic boundary value problem
\begin{equation}\label{eq:semilinear}
\begin{cases}
-\Delta u_\eps + q\, u_\eps^p = 0, & \text{in }  \Omega, \\[6pt]
u_\eps|_{\p\Omega} = \displaystyle \sum_{j=1}^m \eps_j f_j, &  \text{on } \p\Omega
\end{cases}
\end{equation}
where $\eps = (\eps_1,\dots,\eps_m) \in \mathbb{R}^m$ are small parameters and $f_j \in C^{2,\alpha}(\p\Omega)$, with $0<\alpha<1$, are prescribed boundary data. 
Below, the function $q\in C^{\alpha}(\Omega, \mathbb{R})$ is a called the potential function.
For sufficiently small $\eps$ equation \eqref{eq:semilinear} admits a unique small solution $u_\eps \in C^{2,\alpha}(\Omega)$, see \cite{lassas2021inverse}. Moreover, it can be shown that the solution map $S: C^{2,\alpha}(\p\Omega)\supset B_{\eps}(0)\to C^{2,\alpha}(\Omega)$, $S(f)=u_f$, mapping the boundary value $f\in B_{\eps}(0):=\{ g\in C^{2,\alpha}(\p\Omega)\mid \Vert g\Vert_{C^{2,\alpha}(\p\Omega)} \leq \eps\}$ to the corresponding solution of \eqref{eq:semilinear}, is $C^\infty$-map in the Fr\'echet sense in the neighbourhood of the origin in $C^{2,\alpha}(\p\Omega)$. This justifies the differentiation of \eqref{eq:semilinear} with respect to the small parameters $\eps_j$, and in the sequel we will understand any (partial) derivatives of $u_\eps$ with respect to $\eps$ to mean Fr\'echet derivatives of the solution map $S$.

Next, let us denote
\[
v_j := \left. \frac{\p u_\eps}{\p \eps_j}\right|_{\eps=0}.
\]
Then differentiating \eqref{eq:semilinear} with respect to $\eps_j$ and evaluating at $\eps_1=\ldots=\eps_m=0$ yields
\[
-\Delta v_j = 0 \quad \text{in }\Omega, 
\qquad v_j|_{\p\Omega} = f_j,
\]
Thus, each $v_j$ is a harmonic function independent of $q$.
Let then $\sigma=(\sigma_1,\dots,\sigma_m)$ be a multi-index with $|\sigma|=p$ and define
\[
w :=  D_{\eps_1,\ldots,\eps_m}^\sigma u_\eps \big|_{\eps_1=\ldots=\eps_m=0}=
 \frac {\partial^{\sigma_1}} {\partial \eps_1^{\sigma_1}}\dots
 \frac {\partial^{\sigma_m}} {\partial \eps_m^{\sigma_m}}
 u_{\eps_1,\ldots,\eps_m} \big|_{\eps_1=\ldots=\eps_m=0}.
\]
Differentiating equation \eqref{eq:semilinear} $p$ times with respect to $\eps_1,\ldots,\eps_m$ and evaluating at $\eps_1=\ldots=\eps_m=0$, we obtain the linear problem
\begin{equation}\label{eq:w_equation}
\begin{cases}
-\Delta w + p! \, q \displaystyle\prod_{i=1}^m v_i^{\sigma_i} = 0, & \text{in }  \Omega, \\[6pt]
w = 0, & \text{on } \p\Omega.
\end{cases}
\end{equation}
Equation \eqref{eq:w_equation} shows that the $p$-th order linearization $w$ solves a linear elliptic PDE with source term given by $q(x)$ times a product of harmonic functions $v_j$. 

Finally, to obtain a useful integral identity, we let $v_0 \in H^1(\Omega)$ be any harmonic function;
\[
-\Delta v_0 = 0 \quad \text{in }\Omega.
\]
Multiplying \eqref{eq:w_equation} by $v_0$ and integrating by parts yields
\begin{equation}\label{eq:integral_id}
p! \int_\Omega q\, v_0 \prod_{i=1}^m v_i^{\sigma_i} \, \dd x 
= \int_{\p\Omega} f_0 \, \p_\nu w \, \dd S,
\end{equation}
where $f_0:=v_0\big|_{\p\Omega}$ and $\p_\nu$ denotes the outward normal derivative on $\p\Omega$. The right-hand side of \eqref{eq:integral_id} is determined by the Dirichlet-to-Neumann map of the nonlinear problem. Hence,  knowledge of the DN map determines the integrals
\[
\int_\Omega q\, v_0 \prod_{i=1}^m v_i^{\sigma_i} \, \dd x,
\]
for arbitrary choices of boundary data $\{f_j\}$ and harmonic test functions $v_0$. This provides the key identity for recovering the potential $q$ from the boundary measurements.

\section{Numerical solution of the Nonlinear Schr\"odinger equation} \label{section:newton}

We consider the nonlinear boundary value problem
\begin{equation}\label{eq:schrodinger}
\begin{cases}
     -\Delta u + q\, u^p = 0 \quad & \text{in } \Omega, \\ 
  u = f \quad & \text{on } \p \Omega,
\end{cases}
\end{equation}
where $\Omega \subset \mathbb{R}^d$ is a bounded domain, $p\in\mathbb{N}_{\ge 2}$ the power of the nonlinearity, $q\in C^{\alpha}(\Omega, \mathbb{R})$ is a given potential function, and $f\in C^{2,\alpha}(\p\Omega, \mathbb{C})$ prescribes the Dirichlet boundary conditions.  
The nonlinearity in \eqref{eq:schrodinger} makes the problem nontrivial to solve numerically. We adopt a Newton-type approach applied to its variational formulation.
For this purpose, let $V=\{ u\in H^1(\Omega, \mathbb{C})\mid u\big|_{\p\Omega}=f\}$ and
\[
V_0 := H_0^1(\Omega)= \{ v \in H^1(\Omega,\mathbb{C}) \mid v|_{\p\Omega} = 0 \}.
\]
Multiplying \eqref{eq:schrodinger} by a test function $v \in V_0$ and integrating by parts, we obtain the weak formulation: find $u \in V$ such that
\begin{equation}\label{eq:weak_form}
F(u)(v) := \int_\Omega \nabla u \cdot \overline{\nabla v} \, \dd x 
          + \int_\Omega q(x)\, u^p\, \overline{v} \, \dd x = 0, 
\qquad \forall v \in V_0.
\end{equation}
Here $\overline{z}$ denotes the complex conjugate of $z\in\C$. Therefore $F(u)(v)$ is a nonlinear functional in $u$ that vanishes at the solution.
It can be noted that for each fixed $u\in V$ the conjugate linear functional $F(u):V_0\to\C$ is bounded; the proof is included in Appendix~\ref{app: Newton justified}, as the argument is standard.

To solve \eqref{eq:weak_form}, we apply Newton’s method.  
Given an iterate $u_n \in V$, we seek an increment $\delta u \in V_0$ such that
\begin{equation}\label{eq:newton_step}
  DF(u_n)[\delta u; v] = -F(u_n)(v),
  \qquad \forall v \in V_0,
\end{equation}
where $DF(u_n)$ denotes the Fréchet derivative of $F$ with respect to $u$, evaluated at $u_n$.
Computing this derivative, we obtain the Jacobian sesquilinear form
\begin{equation}\label{eq:jacobian}
DF(u_n)[\delta u; v] = J(u_n; \delta u, v) 
  = \int_\Omega \nabla \delta u \cdot \overline{\nabla v} \, \dd x 
  + \int_\Omega p q(x)\, u_n^{p-1} \, \delta u \, \overline{v} \, \dd x.
\end{equation}
It can be verified (see Appendix~\ref{app: Newton justified}) that for a fixed and sufficiently small $u_n\in H^1(\Omega)$, the Jacobian $J$ is a bounded and coercive sesquilinear form on $V_0\times V_0$. By the Lax-Milgram theorem, it is therefore invertible as an operator from $V_0$ to $V_0'$ and the linearized problem \eqref{eq:newton_step} admits a unique solution $\delta u\in V_0$.

Thus, the Newton step consists of solving the linearized problem \eqref{eq:newton_step} with the Jacobian \eqref{eq:jacobian}, and updating
\[
u_{n+1} = u_n + \delta u.
\]
The iterations are repeated until convergence, typically monitored by the residual norm
\[
\| F(u_n) \|_{V_0'} < \eps_\mathrm{tol}
\]
for a prescribed tolerance $\eps_\mathrm{tol} > 0$.  
This procedure ensures quadratic convergence (see \cite{Kelley1995}, \cite{OrtegaRheinboldt1970}), provided the initial guess $u_0$ is sufficiently close to the solution $u$, and the Jacobian remains coercive through the iteration.

\begin{theorem}[Well-posedness of the Newton step and convergence]
\label{thm:newton_convergence}
Let $F:V \to V_0'$ be the nonlinear residual associated with the boundary value
problem \eqref{eq:schrodinger} on a smooth bounded domain $\Omega$, with $V$ and $V_0$
as defined above.
\begin{enumerate}
\item[(i)]  
If $u_n$ is sufficiently close to the exact solution $u$, then the 
linearized problem
\[
DF(u_n)[\delta u] = -F(u_n)
\quad \text{in } V_0'
\]
admits a unique solution $\delta u\in V_0$.

\item[(ii)]  
If the initial guess $u_0\in V$ is sufficiently close to $u$, 
the Newton iteration
\[
u_{n+1}=u_n+\delta u_n,
\qquad 
\delta u_n \text{ solving the linearized step},
\]
is well defined and converges to $u$ in $V$.
\end{enumerate}
\end{theorem}

\paragraph{Implementation notes}

In practice, this method is implemented in a finite element framework. At each Newton iteration, a linear system is assembled from the Jacobian and the residual, and subsequently solved. The nonlinearity enters both the residual and the Jacobian, so care must be taken to update these at each step. In our implementation, the Python-based FEniCX environment handles the automatic differentiation of the variational form to compute the Jacobian (see \cite{dolfinx_NewtonSolver}, \cite{fenicsx_tutorial_newton}).

For the numerical simulations \ref{section: Numerical_examples}, the nonlinear Schrödinger equation is discretized using a conforming finite element method on a circular domain $\Omega = \{ x \in \mathbb{R}^2 : |x| \leq 1 \}$.
The mesh is generated from a quasi-uniform triangulation of the disk obtained from a radial resolution of $64$ points, corresponding to an average mesh size
$h \approx \tfrac{2}{64} \simeq 3.1\times 10^{-2}$. The solution is approximated
in a Lagrange $P_3$ finite element space. The nonlinear system is solved using Newton method with a residual based stopping criterion and a
relative tolerance of $10^{-8}$. For all simulations, Newton's method converges
in a small number of iterations (typically between $2$ and $3$ iterations), compared to a typical Banach fixed-point iteration, which can take 25--35 steps.

\section{Reconstruction of the potential}\label{section:rec method of q}

In this section we present the theoretical framework for reconstructing the potential \(q\) in the nonlinear equation \eqref{eq:nonlinear} with \(p=2\).
In Section~\ref{section:HOL}, we applied a higher order linearization; with a suitable choice of boundary data see Section~\ref{boundary choice} this returns explicitly the Fourier transform \(\mathcal{F}(q)\).
Section~\ref{subsection: E} then develops a discretized formulation of the inverse problem, and Section~\ref{subsection: q_rec} describes its numerical solution via Tikhonov or total-variation (TV) regularization, selected according to the expected smoothness of \(q\).

We consider the nonlinear boundary‐value problem
\begin{equation}\label{eq:nonlin}
\begin{cases}
  -\Delta u + q\,u^2 = 0,
  & \text{in } \Omega,\\
  u = \eps\,f,
  & \text{on } \p\Omega.
\end{cases}
\end{equation}
As in Section~\ref{section:HOL}, let us denote
\[
  v \;=\;\left.\frac{\p u}{\p\eps}\right|_{\eps=0},
  \qquad
  w \;=\;\left.\frac{\p^2 u}{\p\eps^2}\right|_{\eps=0}.
\]

By applying higher order linearization method (see Section~\ref{section:HOL} with $p=2$ in \eqref{eq:integral_id}), we have the integral identity:  
\begin{equation}\label{eq:identity}
  \int_{\p\Omega}f_0\,\p_\nu w\,\dd s
  = 2\int_\Omega q\,v^2\,v_0\,\dd x.
\end{equation}

\begin{remark}
Note that for the purpose of numerical simulations it is possible to use the  identity 
\begin{equation}\label{eq:Lambda2}
\bigl\langle \Lambda_q(f),\,f_0\bigr\rangle_{\partial\Omega}
= \int_\Omega q\,u^2 \;v_0 \,\dd x,
\end{equation}
to get an easier way to evaluate the left hand side of the equation ~\eqref{eq:identity} as
\begin{equation}\label{eq:Widentity}
\int_{\partial\Omega}v_0\,\partial_\nu w\,\dd s
=\frac{\partial^2}{\partial\varepsilon^2} \int_\Omega v_0\,q\,u^2\,\dd x.
\end{equation}
Indeed, in simulations one could use the right-hand side of \eqref{eq:Widentity} as the synthetic data. The benefit here is that one does not need to numerically compute the normal derivative of $u$ on the boundary $\p\Omega$, rather the integrals of $v_0qu^2$ over $\Omega$ suffice as the data.
\end{remark}

\subsection{Choice of the harmonic functions and boundary values} \label{boundary choice}
There are many possible ways to choose the auxiliary function $v_0$ and the boundary value $f$. A first choice might be to select $v_0$ and $f$ so that we have the Calder\'on's exponentials 
\[
  v_0(x) \;=\; \ee^{(-\zeta - \ii\,\xi)\cdot x/2}, 
\qquad
  v(x) \;=\; \ee^{(\zeta - \ii\,\xi)\cdot x/4},
\]
where $|\xi| = |\zeta|$ and $\xi\cdot \zeta  = 0$. Let us denote $f_0:=v_0\big|_{\p\Omega}$. Then from \eqref{eq:Widentity} we obtain
\begin{equation}\label{eq:fourier_q_usual}
\frac{\p^2}{\p\eps^2}\int_{\p\Omega} f_0\,\Lambda_q(\eps f)\,\dd s\Big|_{\eps=0}
\;=\;
2\int_\Omega q(x)\,\ee^{- \ii\,\xi\cdot x}\,\dd x.
\end{equation}
By Fourier inversion, this provides an explicit reconstruction formula for the potential $q$ given measurements of $\Lambda_q$.

Instead of varying $v_0$, we opt to use combinations of boundary values of $f_\pm=(v_1\pm v_2)\big|_{\p\Omega}$, which yields two corresponding solutions $u_+$ and $u_-$.
Then the following algebraic identity holds:
\[
(v_1+v_2)^2 - (v_1-v_2)^2  = 4v_1v_2.
\]
Here, the choice of Calder\'on exponentials
\[
v_1(x)=\frac12 \ee^{(-\zeta -\ii \xi)\cdot x/2}\quad\text{and}\quad 
v_2(x)=\frac12 \ee^{(\zeta -\ii \xi)\cdot x/2}
\]
yields
\[
v_1v_2 =\frac14 \ee^{-\ii \xi\cdot x}.
\]
Therefore, if we can obtain the \emph{two} measurements (for each $\xi$):
\[
\frac{d^2}{d\eps^2}\int_{\p\Omega} f_0 \Lambda_q(\eps f_\pm) \dd s \big|_{\eps=0}
\;
=
2\int_{\Omega} qv_0(v_1\pm v_2)^2 \dd x
\]
their difference produces
\begin{equation}\label{eq:fourier_q_usual_2}
\frac{d^2}{d\eps^2}\int_{\p\Omega} v_0 \Lambda_q(\eps (f_+ -f_-)) \dd s\big|_{\eps=0}
\;
=
2\int_{\Omega} 4qv_0v_1v_2 \dd x
=
2\int_{\Omega}q(x)\ee^{-\ii \xi\cdot x}\dd x,
\end{equation}
where we chose $v_0\equiv 1$ (compare the LHS to \eqref{eq:fourier_q_usual}).

This latter choice of boundary values $f_\pm$ is motivated physically: suppose there is a measurement device, which for each input $f$ simply integrates the boundary value $\partial_\nu u\big|_{\p\Omega}$ and returns a complex number. We \emph{do not need to know the full DN-map} (operator), rather it suffices to measure the functional $\psi: C^{2,\alpha}(\p\Omega)\supset B_\eps(0)\to\C$, mapping
\begin{equation}\label{eq: measurement function}
\psi:f\mapsto \langle \Lambda_q(f),1\rangle_{L^2(\p\Omega)}=\int_{\p\Omega} \Lambda_q(f)dS_x.
\end{equation}
Following \cite{WaveStability}, we call $\psi$ a \emph{measurement function}.

\begin{remark}
Some care is required when making numerical differences of the type
$u_{+}-u_{-}$, to avoid canceling significant digits. This
is already visible at the level of the corresponding Calder\'on exponentials.
Let $v_1 = e^{(\zeta + i\xi)\cdot x/2}$ and $v_2 = e^{(-\zeta + i\xi)\cdot x/2}$,
and denote $s=\zeta\cdot x$. Then
\[
(v_1+v_2)^2 = \ee^{i\xi\cdot x}\,(\ee^{s}+\ee^{-s})^2 = \ee^{i\xi\cdot x}\,4\cosh^2(s),
\]
and
\[
(v_1-v_2)^2 = \ee^{i\xi\cdot x}\,(\ee^{s}-\ee^{-s})^2 = \ee^{i\xi\cdot x}\,4\sinh^2(s).
\]
For $|s|\gg 1$ we have
\[
\cosh^2(s),\ \sinh^2(s)\sim \tfrac14 e^{2|s|},
\]
so each individual quantity $(v_1\pm v_2)^2$ is of magnitude $\sim \ee^{2|s|}$,
while their difference
\[
(v_1+v_2)^2 - (v_1-v_2)^2
= 4 \ee^{i\xi\cdot x}
\]
remains $O(1)$. Thus, the relative size of the difference compared to the
operands is approximately
\[
\frac{|(v_1+v_2)^2 - (v_1-v_2)^2|}
     {\max\{|(v_1+v_2)^2|,\,|(v_1-v_2)^2|\}}
\approx 4 \ee^{-2|s|},
\]
which becomes smaller than machine precision when $|s|$ is too large.
In our numerical experiments, we ensure numerical stability by choosing
$\zeta$ and the computational grid so that $4 \ee^{-2|\zeta\cdot x|} \ge 2.22\times10^{-16}$,
which is satisfied by imposing $|\zeta\cdot x|\le 18.7$ (we use the more comfortable bound $|\zeta\cdot x|\le 10$) on the
computational domain.
\end{remark}

We also mention that purely real boundary values can also be used. The technique is analogous to those above, so we only briefly outline the strategy. Consider
\[
\begin{cases}
  -\Delta u+ q\,u^2 = 0,
  & x\in\Omega,\\
  u = \eps_1 f_1 + \eps_2 f_2,
  & x\in\p\Omega.
\end{cases}
\]
Now it can be seen that the first linearizations $v_j:=D_\eps u\big|_{\eps_1=\eps_2=0}$ in $\Omega$, $j=1,2$, are harmonic with the Dirichlet boundary values $f_j$. The second linearizations $w_{\alpha}=D_{\eps_1,\eps_2}^\alpha u\big|_{\eps_1=\eps_2=0}$ with the multi-index $\alpha\in\N^2$ with $|\alpha|=2$ satisfy $-\Delta w_\alpha+2qv_1^{\alpha_1}v_2^{\alpha_2}=0$ in $\Omega$ with zero Dirichlet boundary values. Therefore, the solution $w$ to $-\Delta w + 2qv^2=0$ in $\Omega$, where $ v$ is harmonic with complex boundary value $f$, can be obtained from real boundary values $f_1=\Re(f)$ and $f_2=\Im(f)$ as the linear combination
\[
w = w_{(2,0)} + 2i w_{(1,1)}- w_{(0,2)}.
\]
This argument holds in general with power-type nonlinearities.

\subsection{Numerical reconstruction of the potential from its Fourier transform} \label{subsection: E}

We embed the computation domain \(\Omega\subset\mathbb{R}^2\) into an axis-aligned rectangle, for \((N_x,N_y)\in\mathbb{N}^2\) denoting the numbers of cells in
\(x\) and \(y\) directions, we define
\[
B=[x_{\min},\,x_{\min}+N_x h]\times [y_{\min},\,y_{\min}+N_y h]
\quad\text{with}\quad \overline{\Omega}\subset B,
\]
and first construct a uniform Cartesian mesh of \(B\) by square cells \(\{p_{m,n}\}\) of side \(h>0\).
Each cell is
\begin{equation}
\begin{split}
p_{m,n}
&=[x_m,x_{m+1}]\times[y_n,y_{n+1}]\\
&= [\,x_{\min}+m h,\; x_{\min}+(m+1)h\,]\times[\,y_{\min}+n h,\; y_{\min}+(n+1)h\,]
\end{split}
\end{equation}
with
\begin{equation*}
(m,n)\in\{0,\dots,N_x-1\}\times\{0,\dots,N_y-1\}. 
\end{equation*}
We then restrict the reconstruction grid to $\Omega$, retaining only the cells whose centers
\[
c_{m,x}:=x_{\min}+\Big(m+\tfrac12\Big)h,\qquad
c_{n,y}:=y_{\min}+\Big(n+\tfrac12\Big)h
\]
satisfy \((c_{m,x},c_{n,y})\in \Omega\). The number of pixels retained is then $N'=N_x'\times N_y'$, in obvious notation. Removing cells outside $\Omega$ reduces the size of the linear system later on and prevents the inversion from reconstructing non-physical mass outside of $\Omega$, which might degrade the reconstruction.

Next, approximate a (complex-valued) \( f \colon \Omega \to \mathbb{C} \)  by a piecewise constant function over each cell using the basis expansion:
\[
\mathbf{f}(X) = \sum_{m,n} f_{m,n} \, \varphi_{m,n}(X),
\]
where the indicator basis functions \( \varphi_{m,n} \) are defined by
\[
\varphi_{m,n}(X) = 
\begin{cases}
1, & \text{if } X \in p_{m,n}, \\
0, & \text{otherwise}.
\end{cases}
\]

To compute the Fourier transform of \( f \), we evaluate the contribution of each cell \( p_{m,n} \) at a discrete set of frequencies \( \xi_k  \in \mathbb{R}^2 \), resulting in the expression
\[
\widehat{f}(\xi_k)=\int_{p_{m,n}} f(X)\, \mathrm{e}^{- \ii \xi_k \cdot X} \, \mathrm{d}X 
= \sum_{m,n}f_{m,n} \int_{p_{m,n}} \mathrm{e}^{- \ii \xi_k \cdot X} \, \mathrm{d}X.
\]
For isotropic resolution, we sample the frequency points \( \xi_k \) in polar coordinates, that is,
\begin{equation}\label{polar}
   \xi_k = (r_\ell \cos\theta_j, \; r_\ell \sin\theta_j) = (\xi_{k,x}, \xi_{k,y}), 
\end{equation}
where \( r_\ell \in [0, R_{\max}] \) denotes the radial coordinate and \( \theta_j \in [0, 2\pi) \) the angular coordinate; here
\(\xi_{k,x}\) and \(\xi_{k,y}\) denote the corresponding Cartesian components. 
We choose a finite number of radii \( N_r \) and angles \( N_\theta \), producing a total of \( N_r \times N_\theta \) sampling points that form a polar grid in the Fourier domain.  

We next build the matrix \( E \in \mathbb{C}^{K \times N'} \), where $N' = N_x' \times N_y'$ and $K = N_r \times N_\theta$, and whose entries \( E_{(k_x, k_y),(m,n)} \) represent the contribution of the cell \( p_{m,n} \) to the Fourier coefficient associated with frequency \( \xi_k \)

If \( \xi_{k,x}  \neq 0 \) and \( \xi_{k,y} \neq 0 \)
\begin{equation*}
\begin{aligned}
E_{(k_x, k_y),(m,n)} &= \int_{x_m}^{x_{m+1}} \int_{y_n}^{y_{n+1}} \mathrm{e}^{- \ii (\xi_{k,x} x + \xi_{k,y} y)} \, \mathrm{d}y \, \mathrm{d}x \\
&= \frac{-1}{ \xi_{k,x} \xi_{k,y}} 
\left( \mathrm{e}^{- \ii \xi_{k,x} h} - 1\right)
\left( \mathrm{e}^{- \ii \xi_{k,y} h} - 1\right)\mathrm{e}^{- \ii (\xi_{k,x} x_{m} + \xi_{k,y} y_{n})},
\end{aligned}
\end{equation*}
where \( x_m = x_{min} + m h \) and \( y_n = y_{min} + n h \), so
\begin{equation*}
\begin{aligned}
E_{(k_x, k_y),(m,n)}=h^2\sinc(\frac12\xi_{k,x}h)\sinc(\frac12\xi_{k,y}h)e^{-i(\xi_{k,x}c_{m,x} + \xi_{k,y}c_{n,y})}
\end{aligned}
\end{equation*}
where the function $\sinc(t)=\frac{\sin(t)}{t}$, with the continuous limit $\sinc(t)=1$ at $t\to 0$ is used in the numerical implementation for stability.

If \( \xi_{k,x} = 0 \) (or \( \xi_{k,y} = 0 \) by symmetry),
\begin{equation*}
\begin{aligned}
E_{(k_x, k_y),(m,n)} =h^2e^{-i\xi_{k,y}c_{n,y}}\sinc(\frac12 \xi_{k,y}h),
\end{aligned}
\end{equation*}

Else,
\begin{equation*}
\begin{aligned}
E_{(k_x, k_y),(m,n)} = h^2 .
\end{aligned}
\end{equation*}

To simplify notation and enable matrix computations, we flatten the two dimensional indices \( (m,n) \) (resp. $(k_x, k_y)$) into a single global index \( I \in \{0, \dots, N' - 1\} \) (resp. \( K \in \{0, \dots, K_r K_\theta - 1\} \)) using row-major ordering:
\[
I = n \cdot N_x + m \quad(\text{resp} \quad K = k_r\cdot K_\theta + k_\theta)
\]

The discrete Fourier transform vector \( F \in \mathbb{C}^K \) can then be approximated as:
\[
F_k \approx \sum_{I=0}^{N'-1} E_{k,I} \, q_I,
\]
or in matrix form,
\begin{equation} \label{eq:à inverser}
  F = E q,  
\end{equation}
where \( q \in \mathbb{C}^{N'} \) contains the piecewise constant values \( q_I = q_{m,n} \), and \( F \in \mathbb{C}^K \) represents the approximated Fourier coefficients of \( q \) at the frequencies \( \xi_k \). 

\begin{remark}
If the potential function $q$ in \eqref{eq:semilinear} is known to be real-valued then, since $\widehat{q}(-\xi)=\overline{\widehat{q}(\xi)}$, it suffices to use frequency components in the upper half-plane $\xi_y\geq 0$ by taking, for example, $\theta_j\in [0,\pi)$ above.
\end{remark}

We finally touch on the choice of the bandwidth $R_{\max}$ of frequencies $\xi$ and the discretization width $h$ of the reconstruction grid of $q$. We choose $N'\geq K$, so that the system $Eq=F$ is under-determined. The usual Nyquist rule suggests that one can expect good resolution reconstructions in fairly coarse grids with $h\sim \pi/R_{\max}$. However, in practise, we choose much finer discretization ($h\ll \pi/R_{\max}$) and suppress the arising small singular values (see Figure~\ref{fig:singular_values}, below) of $E$ via regularization strategies, such as Tikhonov or TV-regularization, using a priori information about the smoothness of the potential $q$.

\subsection{Reconstruction of the potential from noisy measurements} \label{subsection: q_rec}

In practice, measurement data are inevitably contaminated by noise.  
We model this by adding synthetic complex Gaussian noise to the measured Fourier data. The noise level is chosen relative to the 
maximum absolute value of the measurements.

The measurement is
\begin{equation}\label{int+noise}
\int_{\p\Omega} \p_\nu (u_+-u_-)_\eps \dd s \; + \; \eta(\xi)
\end{equation}
where \( \eta(\xi) \) denotes an additive complex Gaussian noise term associated with the Fourier mode \( \xi \in \mathbb{R}^2 \).  
Here, \( \xi \) represents the spatial frequency variable in the Fourier domain, and the random process \( \eta(\xi) \) models independent noise realizations for each frequency component.

When finite-difference schemes are employed to approximate the second derivative with respect to $\varepsilon$ in equation \eqref{eq:fourier_q_usual}, the reconstruction becomes highly unstable in the presence of noise~\cite{Cullum}.  
This instability arises because finite differences amplify measurement errors, leading to poor recovery of the potential $q(x)$.  
To overcome this limitation, we instead employ more robust differentiation strategy, Savitzky-Golay (see appendix~\ref{Savitzky-Golay}) which provide greater noise tolerance while preserving essential features of the signal.

In the numerical reconstructions of Section~\ref{section: Numerical_examples}, we employ the
Savitzky-Golay filter (originally introduced in the chemistry literature for smoothing
noisy spectra \cite{Savitzky1964}) to smooth the data. The method consists of fitting,
on each sliding window of odd length \(2m{+}1\), a local polynomial of degree \(p\) by
least squares and, in our case, evaluating its derivatives at the window center for
numerical differentiation (For further details, see Appendix~\ref{Savitzky-Golay}). We also tested spline-based smoothing, which yielded
competitive results, but we omit those reconstructions here for brevity and because Savitzky-Golay
requires fewer tuning parameters in our setup.

\begin{figure}[h]
  \centering
  \includegraphics[width=0.8\textwidth]{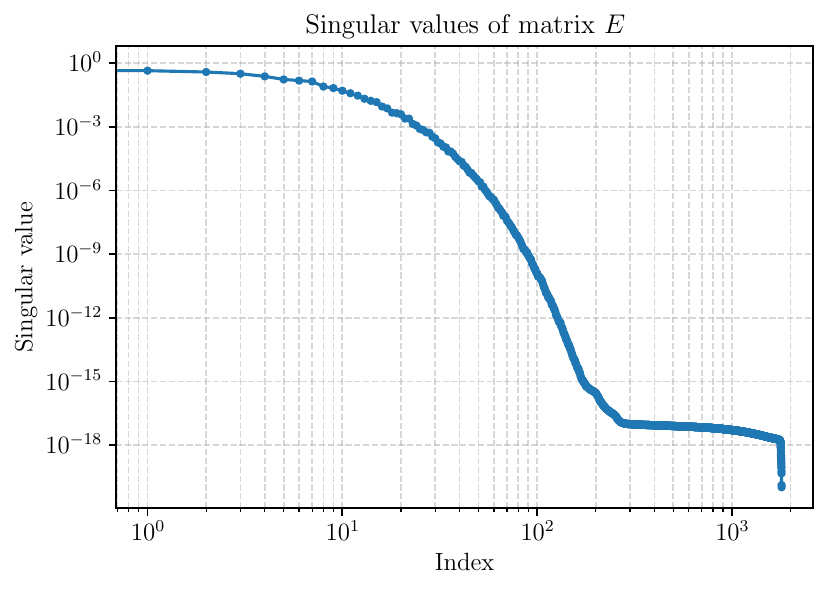}
  \caption{Singular values of the Fourier matrix 
\(E \in \mathbb{C}^{1800 \times 17692}\) in log-log scale 
(\(N_r = 60\), \(N_\theta = 30\), \(N_x = N_y = 150\)).}
  \label{fig:singular_values}
\end{figure}

To recover the potential \( q \) from the discrete linear system~\eqref{eq:à inverser}, we employ regularization techniques that stabilize the inversion and reduce the effects of noise and ill-conditioning. The rapid decay of the singular values (see Fig.~\ref{fig:singular_values}) clearly shows that recovering the potential from Fourier data is a severely ill-posed problem, which justifies the use of Tikhonov or TV regularisation.

Depending on the expected regularity of the underlying potential, we consider two complementary approaches: Tikhonov regularization, suitable for smoothly varying potentials, and Total Variation (TV) regularization, designed for piecewise-smooth or discontinuous potentials. Also, if it is known that recovery of particular frequencies (such as low frequencies) is more stable, this might be used in the regularization strategy. For example, one might use different regularization parameters for low and high frequencies, respectively.

\paragraph{Tikhonov regularization}

To stably recover the unknown vector \( q \in \mathbb{R}^N \) from the possibly ill-posed linear system \( E q = F \), we apply Tikhonov regularization, which introduces a penalization term to enforce stability and incorporate prior knowledge about the solution. The regularized estimate \( q_\lambda \) is obtained as the solution to the minimization problem
\[
q_\lambda = \arg\min_{q \in \mathbb{R}^N} \left\{ \| E q - F \|_2^2 + \lambda \| \Gamma q \|_2^2 \right\},
\]
where \( \lambda > 0 \) is the regularization parameter, and \( \Gamma \in \mathbb{R}^{K \times N} \) is a regularization operator that typically encodes smoothness (e.g., identity or a discrete derivative matrix).
In our setting, we take
\[
\Gamma = \begin{bmatrix}
\mathrm{Id} ,D_x ,D_y
\end{bmatrix}.
\]
where \(D_x\) (resp. \(D_y\)) denotes the discrete derivative operator in the \(x\)-direction (resp. in the \(y\)-direction).

The minimizer \( q_\lambda \) satisfies the normal equations:
\[
\left( E^\top E + \lambda \Gamma^\top \Gamma \right) q_\lambda = E^\top F.
\]
This formulation ensures a stable and well-posed solution even when the original system \( E q = F \) is ill-conditioned.

\paragraph{Total Variation (TV) regularization} 

To promote reconstructions that are piecewise smooth and preserve sharp edges, we employ Total Variation (TV) regularization.  
Given the ill-posed linear system \( E q = F \), the TV-regularized estimate \( q_\lambda \) is defined as the minimizer of the variational problem
\[
q_\lambda = \arg\min_{q \in \mathbb{R}^N}
\left\{ \| E q - F \|_2^2 + \lambda \, \mathrm{TV}(q) \right\},
\]
where \( \lambda > 0 \) controls the trade-off between data fidelity and regularity, and the TV semi-norm is given by
\begin{equation*}
    \mathrm{TV}(q) =
\displaystyle \sum_{i} \left( | \partial_x q_i | + | \partial_y q_i | \right).
\end{equation*}

Unlike the quadratic penalization in Tikhonov regularization, the TV penalty is non-smooth, which allows it to preserve discontinuities in \( q \) (such as edges in images or jumps in piecewise-constant potentials).  
This makes TV particularly effective in inverse problems where the true potential is expected to be sparse in its gradient. TV regularization is implemented using a quadratic optimization approach, following the methodology described in \cite{SiltanenBook}.

\section{Numerical examples}\label{section: Numerical_examples}

In this section, we present several numerical experiments illustrating the reconstruction of the potential $ q(x)$ from noisy boundary measurements, 
based on higher order linearization method described in Section~\ref{section:HOL} that allows the recovery of \(q\) from the 2nd linearization of the measurement function $\psi$ in \eqref{eq: measurement function}, outlined in Sections~\ref{subsection: E} and \ref{subsection: q_rec}. We experimented also on varying $v_0$ and using Calder\'on exponentials $v$. This produces qualitatively similar reconstructions and slightly simpler implementation, but for conciseness these results are not included here.
The experiments are designed to test different smoothness assumptions on \( q \).

\subsection{Computational setup}
All simulations were run on \textit{Lehmus}, the University of Oulu’s (Finland) high performance computing cluster. Each job used 64 CPU cores. The wall clock time is dominated by the
forward solver; once the measurement data are available, the inverse reconstruction is comparatively fast,
particularly with Tikhonov regularization.

The forward problems is implemented using the \texttt{FEniCSx} framework for finite element discretization and nonlinear solvers based on Newton’s method as shown in Section~\ref{section:newton}.  
All simulations are performed in two spatial dimensions on a circular domain (but the choice of domains of other shapes is also possible)
\[
\Omega = \{ (x,y) \in \mathbb{R}^2 : x^2 + y^2 \le 1 \},
\]
with Dirichlet boundary conditions. The synthetic DN-map is simulated by evaluating 
\[
\Lambda_q(f) = (n_x\cdot\nabla u)\big|_{\p\Omega},
\]
where $n_x$ is the outward pointing unit normal vector of $\p\Omega$ at $x$.

In all numerical experiments, the measurement data are corrupted by additive complex Gaussian noise to emulate realistic experimental uncertainty as shown in \eqref{int+noise}.  
The measurement, i.e. the integral from \eqref{int+noise} with a finite parameter $\eps_k$, is denoted by
\[
I(\eps_k) := \int_{\p\Omega} \p_\nu (u_+-u_-)_{\eps_k} \dd s.
\]

We define the noisy data as
\[
I_{\mathrm{noisy}}(\varepsilon_k)
= I(\varepsilon_k)
+ \sigma \big( \eta_{1,k} + i\, \eta_{2,k} \big),
\]
where \( \eta_{1,k}, \eta_{2,k} \sim \mathcal{N}(0,1) \) are independent standard normal random variables, and the noise amplitude is scaled according to
\[
\sigma = \rho \, \max_{\varepsilon_k} | I(\varepsilon_k) |,
\]
with \( \rho > 0 \) representing the prescribed relative noise level.  
In the examples presented below, we fix the noise level to \( \rho = 0.01 \), corresponding to \( 1\% \) relative Gaussian noise, which gives an average signal-to-noise ratio of approximately $ 29.24$ dB.

The second derivative with respect to the perturbation parameter \( \varepsilon \) is approximated using a Savitzky–Golay polynomial filter with \( N_\varepsilon = 64 \) sampling points uniformly distributed in \( \varepsilon \in [-2, 2] \) and a smoothing window of size \(\text{window} = 51\). 

For the Tikhonov method, the regularization parameter is chosen experimentally and fixed to \( \lambda = 10^{-10} \) in Examples~\ref{ex:centre} and~\ref{ex:two}, while for TV method we use \( \lambda = 10^{-7} \) in Example~\ref{ex:ring} and \( \lambda = 10^{-6} \) in Example~\ref{ex:star}. The computational grid consists of \( N_x = N_y = 150 \) nodes in each spatial direction. 

The Fourier data are sampled on a polar grid as explained in \ref{polar}.
In the experiments, we take \( N_r = 60 \) radial samples and \( N_\theta = 30 \) angular directions, with maximum radius \( R_{\max} = 5.0 \), resulting in quite limited bandwidth measurements.

In summary, the reconstruction of $q_{rec}$ reduces to the following regularized problem:
\[
q_{\mathrm{rec}}
=\arg\min_{q}\;
\frac12\Big\|\underbrace{(\mathcal D^{2}_{\mathrm{SG}} I_{\mathrm{noisy}})(0)}_{\text{SG 2nd derivative of }I\text{ at }\varepsilon=0}
-\; E q\Big\|_2^{2}
\;+\;\lambda\,\mathcal R(q),
\]
with the regularizer chosen according to the desired features,
\[
\mathcal R(q)\in\big\{\|q\|_{L^2(\Omega)}^{2}\ \text{(Tikhonov)},\ \mathrm{TV}(q)\ \text{(TV)}\big\}.
\]

\subsection{Test potentials} 

To assess reconstruction accuracy under different regularity regimes, we consider several
synthetic test potentials \(q(x,y)\) ranging from smooth to piecewise constant with jump
discontinuities.
            
\smallskip
\textbf{Example~1. Smooth, localized bump}.\label{ex:centre}

Define the smooth, compactly supported profile
\[
\varphi(x;d)=
\begin{cases}
\exp\!\Big(\dfrac{1}{\big|\tfrac{x}{d}\big|^2-1}+1\Big), & |x|<d,\\[1.2ex]
0, & \text{otherwise},
\end{cases}
\]
where \(d>0\) sets the characteristic width. Using this profile, we build a radially
symmetric potential centered at \((x_0,y_0)\) by
\[
q(x,y;x_0,y_0,d)
:= \varphi\!\Big(\sqrt{(x-x_0)^2+(y-y_0)^2}\,;\,d\Big).
\]

\begin{figure}[H]
  \centering
  \begin{subfigure}{1.0\textwidth}
    \centering
    \includegraphics[width=\textwidth]{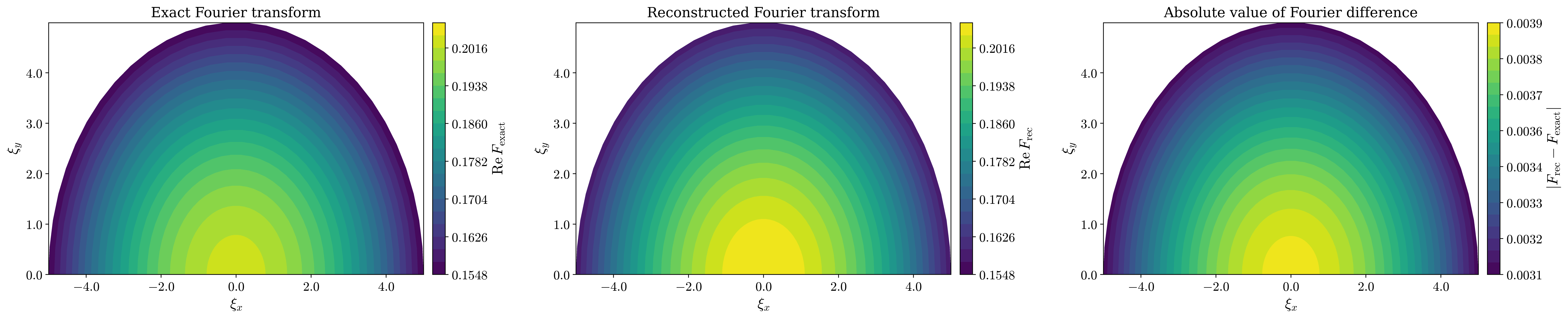}
    \caption{Recovered real Fourier data \(\Re\,\mathcal{F}(q)\) for the centered bump: exact Fourier transform computed via Gauss-Legendre quadrature (left), reconstruction from noisy boundary measurements (middle), and their difference (right).
}
    \label{fig:F_1_bump_centred}
  \end{subfigure}

  \vspace{1em}

  \begin{subfigure}{1.0\textwidth}
    \centering
    \includegraphics[width=\textwidth]{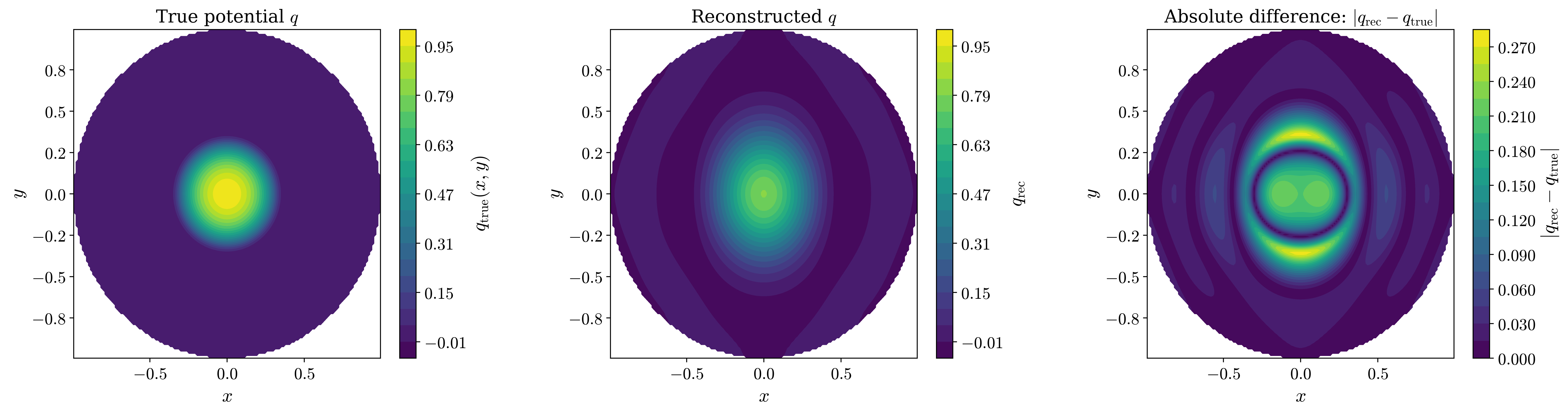}
    \caption{Centered bump \(q(x,y;0,0,0.4)\): ground truth (left), reconstruction (middle), and residual (right). The reconstruction recovers the bump’s location and amplitude well, the $L^2$ error is $0.127$}
    \label{fig:q_1_bump_centred}
  \end{subfigure}
  \caption{Example 1: Smooth, localized bump (centered).}
  \label{fig:ex1_bump}
\end{figure}

Figure~\ref{fig:ex1_bump} reports results for a centered bump at \((0,0)\).
Panel~\ref{fig:F_1_bump_centred} shows the real part of the Fourier transform
of the potential: the reconstruction from boundary data is displayed next to
the exact Fourier data. Panel~\ref{fig:q_1_bump_centred}
presents the spatial reconstruction alongside the ground truth; the location and
shape of the bump are captured accurately.

\begin{figure}[H]
  \centering
  \begin{subfigure}{1.0\textwidth}
    \centering
    \includegraphics[width=\textwidth]{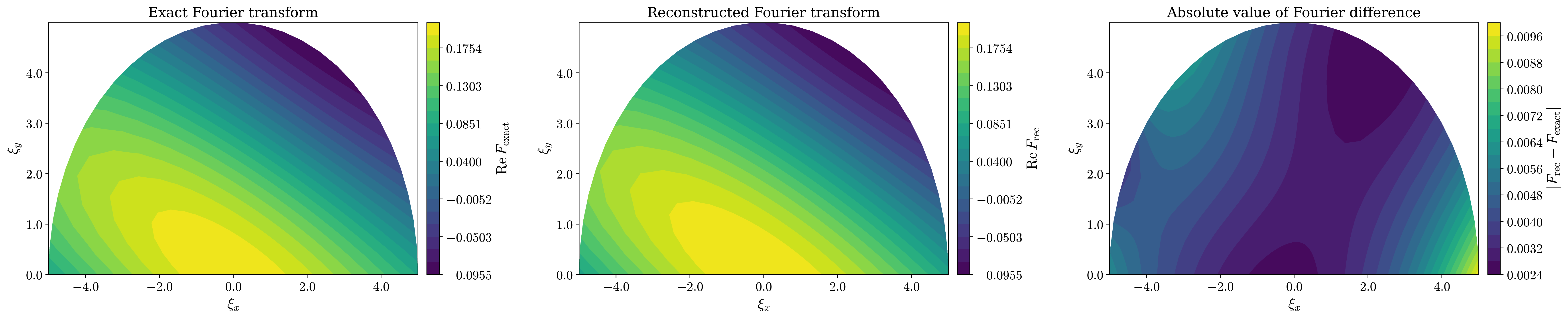}
    \caption{Recovered real Fourier data \(\Re\,\mathcal{F}(q)\) for the off center bump: exact Fourier transform computed via Gauss-Legendre quadrature (left), reconstruction from noisy boundary measurements (middle), and their difference (right).
}
    \label{fig:F_1_bump_off}
  \end{subfigure}

  \vspace{1em}

  \begin{subfigure}{1.0\textwidth}
    \centering
    \includegraphics[width=\textwidth]{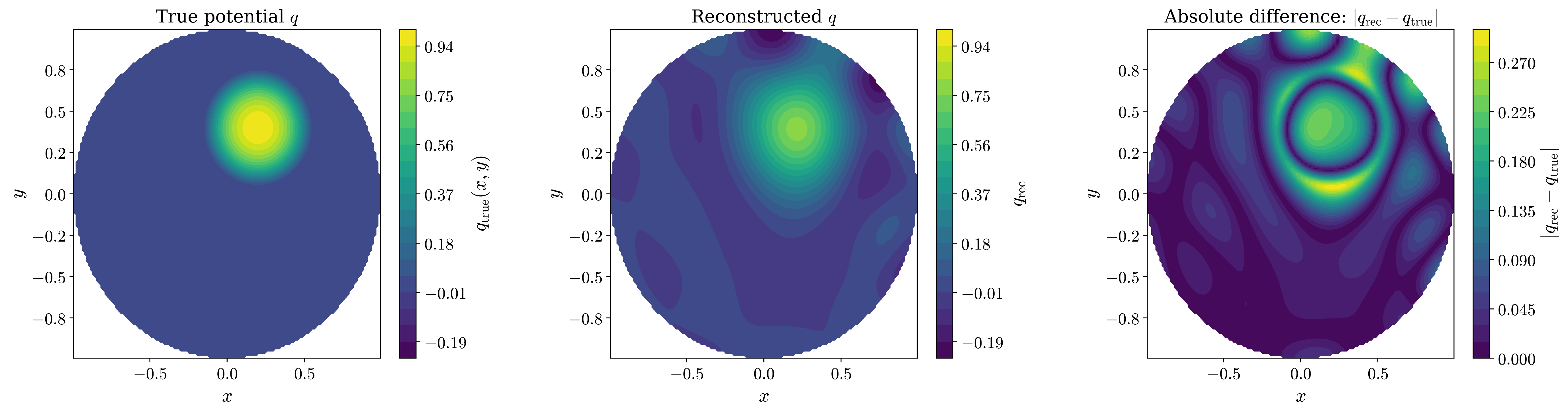}
    \caption{Off center bump \(q(x,y;0.2,0.4,0.4)\): ground truth (left), reconstruction (middle), and residual (right). Both the position and magnitude of the bump are well reproduced ($L^2$ error is $0.147$)}
    \label{fig:q_1_bump_off}
  \end{subfigure}
  \caption{Example 1: Smooth, localized bump (off center).}
  \label{fig:ex1_bump_off}
\end{figure}

Figure~\ref{fig:ex1_bump_off} reports the off–center bump case, located at $(0.2,0.4)$.
Panel~\ref{fig:F_1_bump_off} displays, side by side, the real part of the exact and reconstructed Fourier
transform.
Panel~\ref{fig:q_1_bump_off} shows the reconstructed potential next to the ground
truth; the displacement from the origin and overall shape are recovered effectively.

\smallskip
\textbf{Example~2. Superposition of two bumps.}\label{ex:two}

To introduce mild asymmetry and multiple peaks, we also consider a sum of two localized
bumps:
\[
q(x,y;x_0,y_0,d)
= \varphi\!\Big(\sqrt{(x-x_0)^2+(y-y_0)^2}\,;\,d\Big)
+ \tfrac{1}{2}\,\varphi\!\Big(\sqrt{(x+2x_0)^2+(y+y_0)^2}\,;\,d\Big).
\]

\begin{figure}[H]
  \centering
  \includegraphics[width=1.0\textwidth]{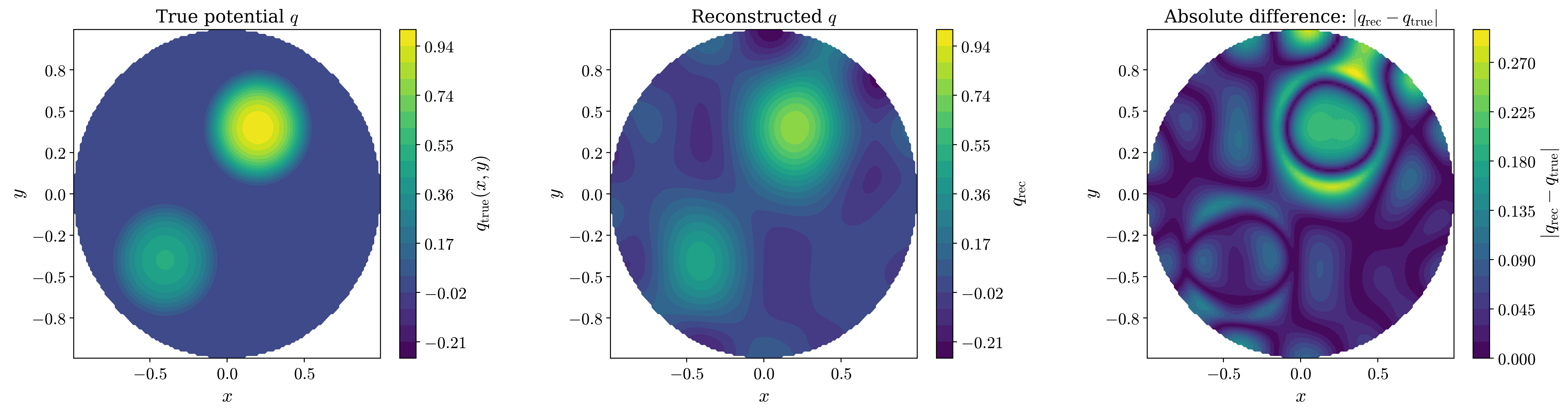}
  \caption{Example 2: Two bump potential \(q(x,y;0.2,0.4, 0.4)\): ground truth (left), reconstruction (middle), and residual (right). The localization and width of the bumps are well reconstructed with Tikhonov regularization ($L^2$ error is $0.153$)}
  \label{fig:ex2_two_bumps}
\end{figure}

The results are shown in Figure~\ref{fig:ex2_two_bumps} and demonstrate a good reconstruction of the two bumps, even with different sizes and amplitudes.

\smallskip
\textbf{Example~3. Discontinuous annulus (ring).}\label{ex:ring}
Let \(\chi_O\) be the indicator of the open annulus
\[
O=\bigl\{(x,y)\in\Omega:\ R_{\mathrm{in}}<\sqrt{x^2+y^2}<R_{\mathrm{out}}\bigr\},
\]
i.e.,
\[
\chi_O(x,y; R_{in},R_{out})=
\begin{cases}
1, & R_{\mathrm{in}}<\sqrt{x^2+y^2}<R_{\mathrm{out}},\\[1mm]
0, & \text{otherwise}.
\end{cases}
\]
and
\[
q(x,y;R_{in},R_{out})=\chi_O(x,y;R_{in},R_{out}).
\]

\begin{figure}[H]
  \centering
  \includegraphics[width=1.0\textwidth]{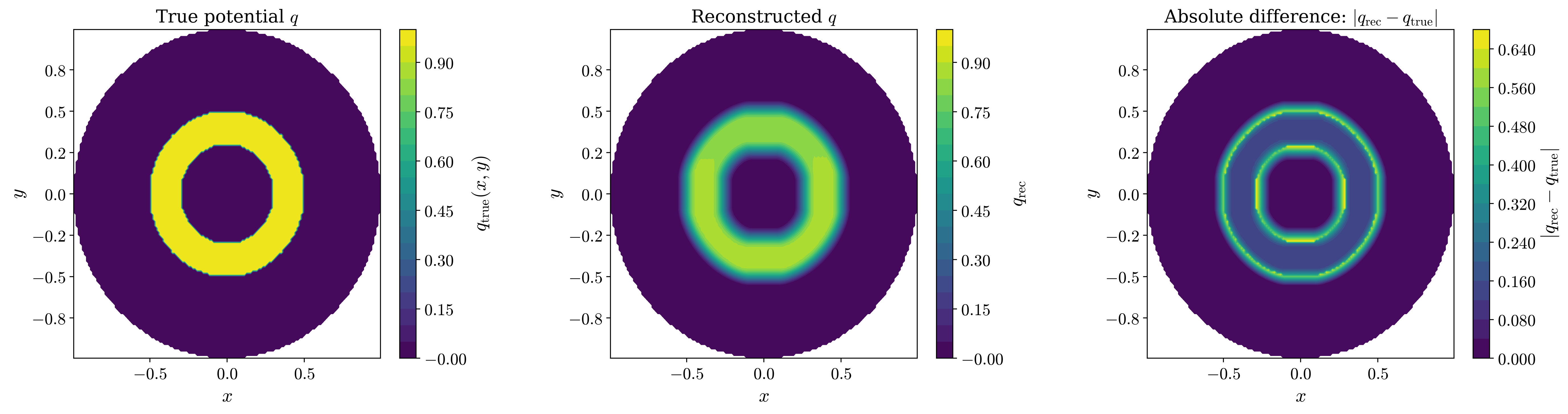}
  \caption{Example 3: Ring potential \(q(x,y;0.3, 0.5)\): ground truth (left), reconstruction (middle), and residual (right). The jump discontinuities are slightly blurred, but piecewise constantness is well-captured with TV regularization ($L^2$ error is $0.229$).}
  \label{fig:ex3_ring}
\end{figure}

Figure~\ref{fig:ex3_ring} shows the results for this example. The reconstruction, obtained with TV regularization, closely captures the sharp inner and outer edges, preserving the jump discontinuities at both radii.

\smallskip
\textbf{Example~4. Discontinuous star shaped set.}\label{ex:star}
Define \(\chi_S\) as the indicator of the star shaped region
\[
S=\Bigl\{(x,y)\in\Omega:\ \sqrt{x^2+y^2}
<R_0+a\cos\big(k\,\theta(x,y)\big)\Bigr\},
\]
where \(\theta(x,y)=\operatorname{atan2}(y,\,x)\), and \(R_0>0\), \(a\in\mathbb{R}\), \(k\in\mathbb{N}\):
\[
\chi_S(x,y; R_0, a, k)=
\begin{cases}
1, & \sqrt{x^2+y^2}<R_0+a\cos\big(k\,\theta(x,y)\big),\\[1mm]
0, & \text{otherwise}.
\end{cases}
\]
and 
\[
q(x,y; R_0, a, k)= \chi_S(x,y; R_0, a, k).
\]

\begin{figure}[H]
  \centering
  \includegraphics[width=1.0\textwidth]{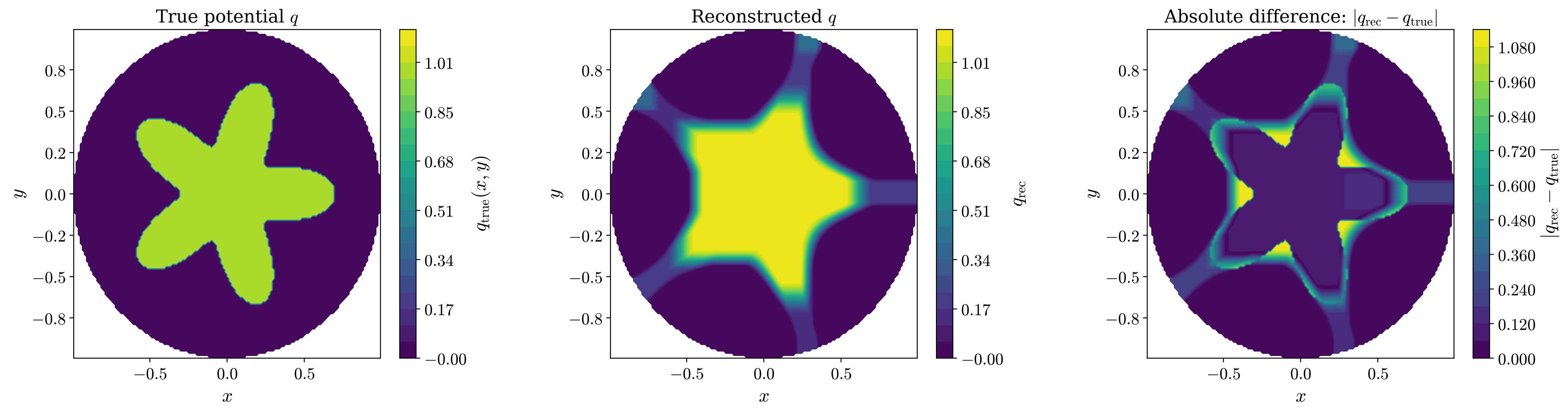}
  \caption{Example 4: Star shaped potential \(q(x,y;0.5,0.2,5)\): ground truth (left), reconstruction (middle), and residual (right). TV regularization recovers the star’s overall shape and sharp edges, but it falls short at the inner corners due to the geometry’s complexity.}
  \label{fig:ex4_star}
\end{figure}

In Figure~\ref{fig:ex4_star}, the reconstruction obtained with TV regularization captures the overall
star shaped contour; however, given the complexity
of this geometry, the reconstruction does not fully recover the inner corners of the star.

\section*{Acknowledgements}
K.E-M. and T.T. were supported by the Finnish Ministry of Education and Culture’s Pilot for Doctoral Programmes (Pilot project Mathematics of Sensing, Imaging and Modelling), the Research Council of Finland (Flagship of Advanced Mathematics for Sensing Imaging and Modelling grant 359186), and Emil Aaltonen
foundation. M.L. was supported by a AdG project 101097198 of the European Research Council, Centre of Excellence of Research Council of Finland, and the
FAME flagship of the Research Council of Finland (grant 359182).
 Views and opinions expressed are 
those of the authors only and do not necessarily reflect those of the European
Union or the other funding organizations.  The European Union  can not be held responsible for them.

\appendix

\section{Justification of the Newton method}\label{app: Newton justified}
To apply the Lax-Milgram theorem, we use standard tools to establish the
boundedness and coercivity of the sesquilinear form \eqref{eq:jacobian}. So, let $\Omega\subset\mathbb{R}^d$ be a bounded Lipschitz domain.

We first show that the linear operator \( F(u_n): V\to \C \) defined in~\eqref{eq:weak_form} is continuous; that is, for any fixed \( u_n \in V \), there exists a constant \( C(n) > 0 \) such that  
\[
|F(u_n)(v)| \leq C(n) \, \|v\|_{H^1(\Omega)} \quad \text{for all } v \in V.
\]
Indeed, by Cauchy-Schwarz inequality and the Sobolev embedding $H^1(\Omega)\subset L^q(\Omega)$ for any $1\leq q<\infty$ (due to boundedness of $\Omega)$ we have
\begin{equation*}
\begin{aligned}
|F(u_n)(v)| 
&\le \|\nabla u_n\|_{L^2(\Omega)} \|\nabla v\|_{L^2(\Omega)}
   + \|q\|_{L^\infty(\Omega)} \|u_n^p\|_{L^2(\Omega)} \|v\|_{L^2(\Omega)} \\
&\le \|\nabla u_n\|_{L^2(\Omega)} \|\nabla v\|_{L^2(\Omega)}
   + \|q\|_{L^\infty(\Omega)} \|u_n\|_{L^{2p}(\Omega)}^p \|v\|_{L^2(\Omega)} \\
&\le \big( \|\nabla u_n\|_{L^2(\Omega)}
   + \|q\|_{L^\infty(\Omega)} C_1 \|u_n\|_{H^1(\Omega)}^p \big)
   \|v\|_{H^1(\Omega)} \\
&\le C(n)\, \|v\|_{H^1(\Omega)},
\end{aligned}
\end{equation*}
with $C(n) = \|\nabla u_n\|_{L^2(\Omega)}
   + \|q\|_{L^\infty(\Omega)} C_1 \|u_n\|_{H^1(\Omega)}^p$.
Hence, \( F(u_n) \) defines a bounded conjugate linear functional $V\to\C$.

To establish the boundedness of \(J(u_n;\cdot,\cdot)\) we, similarly as above, apply Cauchy-Schwarz, generalized H\"older inequality, Poincar\'e inequality, and the Sobolev embeddings \(H^1(\Omega)\hookrightarrow L^{3}(\Omega)\) and \(H^1(\Omega)\hookrightarrow L^{3(p-1)}(\Omega)\):
\begin{equation*}
\begin{aligned}
|J(u_n;u,v)|
&\le \left|\int_\Omega \nabla u \cdot \nabla \overline{v}\,dx\right|
   + \left|\int_\Omega p\,q(x)\,u_n^{p-1}\,u\,\overline{v}\,dx\right| \\[2mm]
&\le \|\nabla u\|_{L^2(\Omega)}\,\|\nabla v\|_{L^2(\Omega)}
   + p\,\|q\|_{L^\infty(\Omega)}\,\|u_n^{p-1}\|_{L^3(\Omega)}\,
      \|u\|_{L^3(\Omega)}\,\|v\|_{L^3(\Omega)} \\[2mm]
&\le \|\nabla u\|_{L^2(\Omega)}\,\|\nabla v\|_{L^2(\Omega)}
   + p\,\|q\|_{L^\infty(\Omega)}\,\|u_n\|_{L^{3(p-1)}(\Omega)}^{p-1}\,
      \|u\|_{L^3(\Omega)}\,\|v\|_{L^3(\Omega)} \\[2mm]
&\le \,\|u\|_{H^1(\Omega)}\,\|v\|_{H^1(\Omega)}
   + C_{0}\,p\,\|q\|_{L^\infty(\Omega)}\,\|u_n\|_{H^1(\Omega)}^{p-1}\,
      \|u\|_{H^1(\Omega)}\,\|v\|_{H^1(\Omega)} \\[2mm]
&\le \Big(1 + C_{0}\,p\,\|q\|_{L^\infty(\Omega)}\,\|u_n\|_{H^1(\Omega)}^{p-1}\Big)\,
   \|u\|_{H^1(\Omega)}\,\|v\|_{H^1(\Omega)} .
\end{aligned}
\end{equation*}
Here  \(C_{0}>0\) is a generic constant depending only on \(\Omega\) and the Sobolev embeddings
\(H^1(\Omega)\hookrightarrow L^{3}(\Omega)\) and \(H^1(\Omega)\hookrightarrow L^{3(p-1)}(\Omega)\)
used in the estimate.

Finally, to prove the coercivity of $J(u_n;\cdot,\cdot)$ we estimate using the Sobolev embedding and Poincar\'e inequality
\begin{equation*} 
    \begin{aligned}
        |J(u_n; u,u)| &\geq |\int_\Omega \nabla u \cdot \nabla \overline{u} \, dx 
  |- |\int_\Omega p q(x)\, u_n^{p-1} \,  u \, \overline{u} \, dx|\\
  & \geq \Vert\nabla u\Vert^2_{L^2(\Omega)} - p \Vert q\Vert_{L^{\infty}(\Omega)} \Vert u_n^{p-1}\Vert_{L^2(\Omega)} \Vert u^2\Vert_{L^2(\Omega)} \\
  & = \Vert\nabla u\Vert^2_{L^2(\Omega)} - p \Vert q\Vert_{L^{\infty}} \Vert u_n\Vert^{p-1}_{L^{2(p-1)}(\Omega)} \Vert u\Vert^2_{L^4(\Omega)} \\
  & \geq C_2  \Vert u\Vert_{H^1(\Omega)}^2 - C_3p \Vert q\Vert_{L^{\infty}(\Omega)} \Vert u_n\Vert^{p-1}_{H^1(\Omega)} \Vert u\Vert^2_{H^1(\Omega)}\\
  & = C_2(1 - \frac{C_3}{C_2}p \Vert q\Vert_{L^{\infty}(\Omega)} \Vert u_n\Vert^{p-1}_{H^1(\Omega)})\Vert u\Vert^2_{H^1(\Omega)}\\
  & \geq \beta_n \|u\Vert^2_{H^1(\Omega)}
    \end{aligned}
\end{equation*}
Above we used a single constant \(C_3\) to represent the Sobolev embedding constants controlling \(\|u_n\|_{L^{2(p-1)}(\Omega)}\) and \(\|u\|_{L^{4}(\Omega)}\) by \(\|\cdot\|_{H^1(\Omega)}\), while \(C_2\) denotes the constant arising from Poincar\'e's inequality.

Note $\beta_n>0$ if
\[
\Vert u_n\Vert_{H^1(\Omega)}^{p-1}< \frac{C_2}{C_3p\Vert q \Vert_{L^\infty(\Omega)}}.
\]
Thus $J(u_n;\cdot,\cdot)$ is bounded coercive sesquilinear form if $u_n\in V$ is sufficiently small. Since $F(u)$ is linear and bounded, the Lax-Milgram theorem applies and the Newton iteration is justified.

\smallskip

\begin{remark}
We recall the standard Newton-Kantorovich convergence result (see Chapter~12, Theorem~12.6.2 of~\cite{OrtegaRheinboldt1970}).  
We have that the Fréchet derivative $DF(u)=J(u)$ satisfies:
\begin{itemize}
    \item[1.] Bound on the inverse,
   \[
     \|J(u_n)^{-1}\|_{\mathcal{L}(V^*,V)} \le M,
   \]
   this follows from Lax-Milgram theorem that implies that \(J(u_n)\) is
  invertible and that 
  \(\|J(u_n)^{-1}\|_{\mathcal{L}(V^*,V)}\le \beta_n^{-1}\), provided that 
    \[
    \Vert u_n\Vert_{H^1(\Omega)}^{p-1}\leq \frac{C_2}{C_3p\Vert q \Vert_{L^\infty(\Omega)}+1}.
    \]
   \item[2.] Lipschitz condition in a neighbourhood of the solution $u$,
   \[
     \|J(u_n)-J(u)\|_{\mathcal{L}(V,V^*)} \le L \|u_n-u\|_{H^1}.
   \]
   This property follows from the fact that the nonlinearity
  \(w \mapsto w^{p-1}\) is Lipschitz on bounded sets together with the
  boundedness of \(q\) and the Sobolev embeddings in
  two dimensions.
\end{itemize}

Let $u$ be the exact solution of $F(u)=0$ and let $(u_n)$ be the Newton iterates
\[
    u_{n+1} = u_n - J(u_n)^{-1} F(u_n).
\]
A Taylor expansion with integral remainder yields
\[
    F(u_n)
    = J(u)e_n + R_n,
    \qquad e_n := u_n - u,
\]
with
\[
    \|R_n\|_{V^*} \le \frac{L}{2}\|e_n\|_V^2.
\]
Using the Newton update and the above properties of $J$, one obtains the classical estimate
\[
    \|e_{n+1}\|_V
    \le M L \|e_n\|_V^2 + \frac{M L}{2}\|e_n\|_V^2
    = \frac{3ML}{2} \|e_n\|_V^2,
\]
which shows quadratic convergence as long as the iterates remain in the neighbourhood where the assumptions are valid.

\end{remark}

\section{Savitzky-Golay regularized higher order differentiation as local polynomial regression} \label{Savitzky-Golay}

Let $\{(x_i,y_i)\}_{i=1}^N$ be samples of a (possibly complex valued) function
$f:\mathbb{R}\to\mathbb{C}$ observed with additive noise:

\begin{equation}
y_i \;=\; f(x_i) + \varepsilon_i,\qquad \mathbb{E}[\varepsilon_i]=0,\quad
\operatorname{Var}(\varepsilon_i)=\sigma^2.
\end{equation}
Assume a uniform grid $x_i = x_0 + i\,\Delta x$. The Savitzky-Golay (SG) method \cite{Savitzky1964,Press1992,Schafer2011} estimates
the $r$-th derivative $f^{(r)}(x_i)$ by fitting, in a small moving window of $2m+1$ points around $x_i$,
a polynomial of degree $p\ge r$ via least squares, and then differentiating the fitted
polynomial analytically at the window center. This is equivalent to applying a
finite impulse response (FIR) convolution filter whose coefficients depend only on
$(p,r,m)$ and the spacing $\Delta x$.

Fix an index $i$ and an odd window of $2m+1$ points,
\[
\mathcal{I}_i=\{i-m,\ldots,i-1,i,i+1,\ldots,i+m\}.
\]

Introduce centered coordinates $t_j = x_{i+j}-x_i = j\,\Delta x$ for $j=-m,\ldots,m$,
and approximate $f$ on the window by a degree-$p$ polynomial $g$ of the form
\begin{equation}
g(t) \;=\; \sum_{k=0}^{p} \beta_k \, t^k.
\end{equation}

The coefficients $\beta_0,\ldots,\beta_p$ are determined by a
least-squares fit to the noisy data $\{y_{i+j}\}_{j=-m}^{m}$, by solving
\[
(\beta_0,\ldots,\beta_p)
   = \arg\min_{\gamma_0,\ldots,\gamma_p}
     \sum_{j=-m}^{m} \bigl|\, \gamma_0 + \gamma_1 t_j + \cdots + \gamma_p t_j^{\,p}
                           - y_{i+j} \bigr|^2 .
\]

The value of the fitted polynomial at the window center is $g(0)=\beta_0$,
and its $r$-th derivative at $t=0$ is given analytically by
\[
g^{(r)}(0) = r!\,\beta_r.
\]
Thus the SG estimate of $f^{(r)}(x_i)$ is obtained by locally replacing
the noisy data with its best-fitting degree-$p$ polynomial and evaluating
its analytic derivative at the center of the window.

\newpage
% ----- Bibliography -----
\bibliography{biblio}
\bibliographystyle{abbrv}

\noindent{\footnotesize E-mail addresses:\\
Khaoula El Maddah: khaoula.elmaddah@oulu.fi (corresponding author)\\
Matti Lassas: matti.lassas@helsinki.fi\\
Teemu Tyni: teemu.tyni@oulu.fi
}

\end{document}